
\input amstex



\def\next{AMS-J}

\def\next{AMSPPT}
\ifx\styname\next \else\input amsppt.sty \relax\fi

\catcode`\@=11

\indenti=3pc

\font@\fourteenbf=cmbx10 scaled \magstep2

\newtoks\sevenpoint@
\def\sevenpoint{\normalbaselineskip9\p@
 \textonlyfont@\rm\sevenrm \textonlyfont@\it\sevenit
 \textonlyfont@\sl\sevensl \textonlyfont@\bf\sevenbf
 \textonlyfont@\smc\sevensmc \textonlyfont@\tt\seventt
  \textfont\z@\sevenrm \scriptfont\z@\sixrm
       \scriptscriptfont\z@\fiverm
  \textfont\@ne\seveni \scriptfont\@ne\sixi
       \scriptscriptfont\@ne\fivei
  \textfont\tw@\sevensy \scriptfont\tw@\sixsy
       \scriptscriptfont\tw@\fivesy
  \textfont\thr@@\sevenex \scriptfont\thr@@\sevenex
   \scriptscriptfont\thr@@\sevenex
  \textfont\itfam\sevenit \scriptfont\itfam\sevenit
   \scriptscriptfont\itfam\sevenit
  \textfont\bffam\sevenbf \scriptfont\bffam\sixbf
   \scriptscriptfont\bffam\fivebf
  \textfont\msbfam\sevenmsb \scriptfont\msbfam\fivemsb
   \scriptscriptfont\msbfam\fivemsb
 \setbox\strutbox\hbox{\vrule height7\p@ depth3\p@ width\z@}%
 \setbox\strutbox@\hbox{\raise.5\normallineskiplimit\vbox{%
   \kern-\normallineskiplimit\copy\strutbox}}%
 \setbox\z@\vbox{\hbox{$($}\kern\z@}\bigsize@1.2\ht\z@
 \normalbaselines\sevenrm\dotsspace@1.5mu\ex@.2326ex\jot3\ex@
 \the\sevenpoint@}

\def\headlinefont@{\sevenpoint}
\def\foliofont@{\sevenrm}


%
%

\define\issueinfo#1#2#3#4{%
  \def\issuevol@{#1}\def\issueno@{#2}%
  \def\issuemonth@{#3}\def\issueyear@{#4}}

\define\originfo#1#2#3#4{\def\origvol@{#1}\def\origno@{#2}%
  \def\origmonth@{#3}\def\origyear@{#4}}

\define\copyrightinfo#1#2{\def\copyrightyear{#1}\def\crholder@{#2}}

\define\pagespan#1#2{\pageno=#1\def\start@page{#1}\def\end@page{#2}}

\define\PII#1{\def\@PII{#1}}

\issueinfo{00}{0}{}{1997}
\originfo{00}{0}{}{1997}
\copyrightinfo{\issueyear@}{American Mathematical Society}
\pagespan{000}{000}
\pageno=1 
\PII{\issn@(\issueyear@)0000-0}


\def\nojourlogo{\let\jourlogo\empty@}

\def\journame{AMS JOURNAL STYLE}
\def\volinfo{Volume \issuevol@, Number \number\issueno@,
  \issuemonth@\ \issueyear@}
\let\pageinfo\empty@
\let\jourlogoextra@\empty@
\let\jourlogoright@\empty@

\def\jourlogofont@{\sixrm\baselineskip7\p@\relax}
\def\jourlogo{%
  \vbox to\headlineheight{%
    \parshape\z@ \leftskip\z@ \rightskip\z@
    \parfillskip\z@ plus1fil\relax
    \jourlogofont@ \frenchspacing
    \line{\vtop{\parindent\z@ \hsize=.5\hsize
      \journame\newline\volinfo\pageinfo\jourlogoextra@\endgraf\vss}%
      \hfil
      \jourlogoright@
    }%
    \vss}%
}

\def\issn#1{\gdef\issn@{#1}}
\issn{0000-0000}


\def\copyrightline@{%
  \rightline{\sixrm \textfont2=\sixsy \copyright\copyrightyear\ \crholder@}}

\def\logo@{\copyrightline@}

\def\keyboarder#1{}%

\let\thedate@\empty@
\def\date#1\enddate{\gdef\thedate@{%
  \eightpoint Received by the editors \ignorespaces#1\unskip.}}

\let\thededicatory@\empty@
\def\dedicatory #1\enddedicatory{\gdef\thededicatory@{{\vskip10pt
  \eightpoint\it \raggedcenter@#1\endgraf}}}

\def\preabstract{}

\let\thecommby@=\empty@
\def\commby #1\endcommby{\gdef\thecommby@{{\vskip1pc
  \eightpoint \raggedcenter@(Communicated by #1)\endgraf}}}

\def\title{\let\savedef@\title
  \def\title##1\endtitle{\let\title\savedef@
    \def\thetitle@{##1}%
    \global\setbox\titlebox@\vtop{\tenpoint\bf
      \raggedcenter@
      \frills@\uppercasetext@{##1}\endgraf}%
    \ifmonograph@ \edef\next{\the\leftheadtoks}%
      \ifx\next\empty@
      \leftheadtext{##1}\fi
    \fi
    \edef\next{\the\rightheadtoks}\ifx\next\empty@ \rightheadtext{##1}\fi
  }%
  \nofrillscheck\title}

\def\author#1\endauthor{%
  \def\theauthor@{#1}%
  \global\setbox\authorbox@
  \vbox{\eightpoint\raggedcenter@
    \uppercasetext@{#1}\endgraf}\relaxnext@
  \edef\next{\the\leftheadtoks}\ifx\next\empty@\leftheadtext{#1}\fi}

\outer\def\endtopmatter{\add@missing\endabstract
  \edef\next{\the\leftheadtoks}%
  \ifx\next\empty@
    \expandafter\leftheadtext\expandafter{\the\rightheadtoks}\fi
  \ifx\thedate@\empty@\else \makefootnote@{}{\thedate@}\fi
  \ifx\thesubjclass@\empty@\else \makefootnote@{}{\thesubjclass@}\fi
  \ifx\thekeywords@\empty@\else \makefootnote@{}{\thekeywords@}\fi
  \ifx\thethanks@\empty@\else \makefootnote@{}{\thethanks@}\fi
  \inslogo@
  \pretitle
  \begingroup 
    \topskip42pt
    \box\titlebox@
  \endgroup
  \preauthor
  \ifvoid\authorbox@\else \vskip2pc\unvbox\authorbox@\fi
  \ifx\thecommby@\empty@\else \thecommby@\fi
  \ifx\thededicatory@\empty@\else \thededicatory@\fi
  \setabstract@
  \ifvoid\tocbox@\else\vskip1.5pcplus.5pc\unvbox\tocbox@\fi
  \prepaper
  \vskip2pcplus1pc\relax
}

\def\setabstract@{%
  \preabstract
  \ifvoid\abstractbox@\else
       \vskip17ptplus.5pc\unvbox\abstractbox@ \fi
}

\def\specialheadfont@{\tenpoint\rm}
\let\specialhead\relax
\outer\def\specialhead#1\endspecialhead{%
  \add@missing\endroster \add@missing\enddefinition
  \add@missing\enddemo \add@missing\endexample
  \add@missing\endproclaim
  \penaltyandskip@{-200}\aboveheadskip
  \begingroup\interlinepenalty\@M\rightskip\z@ plus\hsize
  \let\\\linebreak
  \specialheadfont@\raggedcenter@\noindent
	\uppercase{#1}\endgraf\endgroup\nobreak\vskip\belowheadskip}

\let\subsubhead\relax
\outer\def\subsubhead{%
  \add@missing\endroster \add@missing\enddefinition
  \add@missing\enddemo
  \add@missing\endexample \add@missing\endproclaim
  \let\savedef@\subsubhead \let\subsubhead\relax
  \def\subsubhead##1\endsubsubhead{\restoredef@\subsubhead
    \penaltyandskip@{-50}\subsubheadskip
    {\def\usualspace{\/{\subsubheadfont@\enspace}}%
     \varindent@\subsubheadfont@##1\unskip\frills@{.\enspace}}\ignorespaces}%
  \nofrillscheck\subsubhead}

\long\def\ext #1\endext{\block #1\endblock}

\def\qed{\ifhmode\unskip\nobreak\fi\hfill
  \ifmmode\square\else$\m@th\square$\fi}

\newskip\prexcaskip  \prexcaskip=\medskipamount
\def\xcaheadfont@{\bf}
\outer\def\xca{\let\savedef@\xca \let\xca\relax
  \add@missing\endproclaim \add@missing\endroster
  \add@missing\endxca \envir@stack\endxca
   \def\xca##1{\restoredef@\xca
     \penaltyandskip@{-100}\prexcaskip
        \bgroup{\def\usualspace{{\xcaheadfont@\enspace}}%
        \varindent@\xcaheadfont@\ignorespaces##1\unskip
        \frills@{.\xcaheadfont@\enspace}}%
        \ignorespaces}%
  \nofrillscheck\xca}
\def\endxca{\egroup\revert@envir\endxca
  \par\medskip}

\newdimen\rosteritemsep
\rosteritemsep=.5pc

\newdimen\rosteritemitemwd
\newdimen\rosteritemitemitemwd

\newbox\setwdbox
\setbox\setwdbox\hbox{0.}\rosteritemwd=\wd\setwdbox
\setbox\setwdbox\hbox{0.\hskip.5pc(c)}\rosteritemitemwd=\wd\setwdbox
\setbox\setwdbox\hbox{0.\hskip.5pc(c)\hskip.5pc(ii)}%
  \rosteritemitemitemwd=\wd\setwdbox

\let\rosterskip@\empty@
\def\roster{%
  \envir@stack\endroster
  \edef\leftskip@{\leftskip\the\leftskip}%
  \relaxnext@
  \rostercount@\z@
  \def\item{\FN@\rosteritem@}%
  \def\itemitem{\FN@\rosteritemitem@}%
  \def\itemitemitem{\FN@\rosteritemitemitem@}%
  \DN@{\ifx\next\runinitem\let\next@\nextii@
    \else\let\next@\nextiii@
    \fi\next@}%
  \DNii@\runinitem
    {\unskip
     \DN@{\ifx\next[\let\next@\nextii@
       \else\ifx\next"\let\next@\nextiii@\else\let\next@\nextiv@\fi
       \fi\next@}%
     \DNii@[####1]{\rostercount@####1\relax
       \therosteritem{\number\rostercount@}~\ignorespaces}%
     \def\nextiii@"####1"{{\rm####1}~\ignorespaces}%
     \def\nextiv@{\therosteritem1\rostercount@\@ne~}%
     \par@\firstitem@false
     \FN@\next@}
  \def\nextiii@{\par\par@
    \penalty\@m\rosterskip@\vskip-\parskip
    \firstitem@true}%
  \FN@\next@}

\def\rosteritem@{\iffirstitem@\firstitem@false
  \else\par\vskip-\parskip
  \fi
  \leftskip\rosteritemwd \advance\leftskip\normalparindent
  \advance\leftskip.5pc \noindent
  \DNii@[##1]{\rostercount@##1\relax\itembox@}%
  \def\nextiii@"##1"{\def\therosteritem@{\rm##1}\itembox@}%
  \def\nextiv@{\advance\rostercount@\@ne\itembox@}%
  \def\therosteritem@{\therosteritem{\number\rostercount@}}%
  \ifx\next[\let\next@\nextii@
  \else\ifx\next"\let\next@\nextiii@\else\let\next@\nextiv@\fi
  \fi\next@}

\def\itembox@{\llap{\hbox to\rosteritemwd{\hss
  \kern\z@ 
  \therosteritem@}\hskip.5pc}\ignorespaces}

\def\therosteritem#1{\rom{\ignorespaces#1.\unskip}}

\def\rosteritemitem@{\iffirstitem@\firstitem@false
  \else\par\vskip-\parskip
  \fi
  \leftskip\rosteritemitemwd \advance\leftskip\normalparindent
  \advance\leftskip.5pc \noindent
  \DNii@[##1]{\rostercount@##1\relax\itemitembox@}%
  \def\nextiii@"##1"{\def\therosteritemitem@{\rm##1}\itemitembox@}%
  \def\nextiv@{\advance\rostercount@\@ne\itemitembox@}%
  \def\therosteritemitem@{\therosteritemitem{\number\rostercount@}}%
  \ifx\next[\let\next@\nextii@
  \else\ifx\next"\let\next@\nextiii@\else\let\next@\nextiv@\fi
  \fi\next@}

\def\itemitembox@{\llap{\hbox to\rosteritemitemwd{\hss
  \kern\z@ 
  \therosteritemitem@}\hskip.5pc}\ignorespaces}

\def\therosteritemitem#1{\rom{(\ignorespaces#1\unskip)}}

\def\rosteritemitemitem@{\iffirstitem@\firstitem@false
  \else\par\vskip-\parskip
  \fi
  \leftskip\rosteritemitemitemwd \advance\leftskip\normalparindent
  \advance\leftskip.5pc \noindent
  \DNii@[##1]{\rostercount@##1\relax\itemitemitembox@}%
  \def\nextiii@"##1"{\def\therosteritemitemitem@{\rm##1}\itemitemitembox@}%
  \def\nextiv@{\advance\rostercount@\@ne\itemitemitembox@}%
  \def\therosteritemitemitem@{\therosteritemitemitem{\number\rostercount@}}%
  \ifx\next[\let\next@\nextii@
  \else\ifx\next"\let\next@\nextiii@\else\let\next@\nextiv@\fi
  \fi\next@}

\def\itemitemitembox@{\llap{\hbox to\rosteritemitemitemwd{\hss
  \kern\z@ 
  \therosteritemitemitem@}\hskip.5pc}\ignorespaces}

\def\therosteritemitemitem#1{\rom{(\ignorespaces#1\unskip)}}

\def\endroster{\relaxnext@
  \revert@envir\endroster 
  \par\leftskip@
  \penalty-50 
  \DN@{\ifx\next\Runinitem\let\next@\relax
    \else\nextRunin@false\let\item\plainitem@
      \ifx\next\par
        \DN@\par{\everypar\expandafter{\the\everypartoks@}}%
      \else
        \DN@{\noindent\everypar\expandafter{\the\everypartoks@}}%
      \fi
    \fi\next@}%
  \FN@\next@}

\widestnumber\key{99}

\catcode`\@=13



\catcode`\@=11

\catcode`\@=13

\catcode`\@=11

\def\journame{JOURNAL OF THE\newline
	AMERICAN MATHEMATICAL SOCIETY}

\issn{0894-0347}


\outer\def\endtopmatter{\add@missing\endabstract
  \edef\next{\the\leftheadtoks}%
  \ifx\next\empty@
    \expandafter\leftheadtext\expandafter{\the\rightheadtoks}
  \fi
  \ifx\thedate@\empty@\else \makefootnote@{}{\thedate@}\fi
  \ifx\thesubjclass@\empty@\else \makefootnote@{}{\thesubjclass@}\fi
  \ifx\thekeywords@\empty@\else \makefootnote@{}{\thekeywords@}\fi
  \ifx\thethanks@\empty@\else \makefootnote@{}{\thethanks@}\fi
  \inslogo@
  \pretitle
  \begingroup 
  \topskip44pt
  \box\titlebox@
  \endgroup
  \preauthor
  \ifvoid\authorbox@\else \vskip2pc\unvbox\authorbox@\fi
  \ifx\thededicatory@\empty@\else \thededicatory@\fi
  \ifvoid\tocbox@\else\vskip1.5pcplus.5pc\unvbox\tocbox@\fi
  \prepaper
  \vskip2pcplus1pc\relax
}

\def\enddocument@text{%
  \ifmonograph@ 
  \else
    \nobreak
    \preabstract
    \ifvoid\abstractbox@\else
      \vskip1.5pcplus.5pc\unvbox\abstractbox@
    \fi
    \count@\z@
    \loop\ifnum\count@<\addresscount@\advance\count@\@ne
      \csname address\number\count@\endcsname
      \csname email\number\count@\endcsname
      \csname urladdr\number\count@\endcsname
    \repeat
  \fi
}

\catcode`\@=13

\NoBlackBoxes

\frenchspacing 

\issueinfo{00}
  {0}
  {Xxxx}
  {1997}



\catcode`\@=12

\def\qed{\quad\raise -2pt\hbox{\vrule\vbox to 10pt{\hrule width 4pt
   \vfill\hrule}\vrule}}

\def\cqfd{\penalty 500 \hbox{\qed}\par\smallskip}

\def\rem#1|{\par\medskip{{\it #1}\pointir}}

\def\vspace[#1]{\noalign{\vskip#1}}

\def\Grille{\setbox13=\vbox to 5mm{\hrule width 110mm\vfill}
\setbox13=\vbox{\offinterlineskip
   \copy13\copy13\copy13\copy13\copy13\copy13\copy13\copy13
   \copy13\copy13\copy13\copy13\box13\hrule width 110mm}
\setbox14=\hbox to 5mm{\vrule height 65mm\hfill}
\setbox14=\hbox{\copy14\copy14\copy14\copy14\copy14\copy14
   \copy14\copy14\copy14\copy14\copy14\copy14\copy14\copy14
   \copy14\copy14\copy14\copy14\copy14\copy14\copy14\copy14\box14}
\ht14=0pt\dp14=0pt\wd14=0pt
\setbox13=\vbox to 0pt{\vss\box13\offinterlineskip\box14}
\wd13=0pt\box13}


\def\fleche(#1,#2)\dir(#3,#4)\long#5{%
\noalign{\nointerlineskip\leftput(#1,#2){\vector(#3,#4){#5}}\nointerlineskip
}}


\def\hfl#1#2#3{\smash{\mathop{\hbox to#3{\rightarrowfill}}\limits
^{\scriptstyle#1}_{\scriptstyle#2}}}

\def\gfl#1#2#3{\smash{\mathop{\hbox to#3{\leftarrowfill}}\limits
^{\scriptstyle#1}_{\scriptstyle#2}}}


 \message{`lline' & `vector' macros from LaTeX}
 \catcode`@=11
\def\{{\relax\ifmmode\lbrace\else$\lbrace$\fi}
\def\}{\relax\ifmmode\rbrace\else$\rbrace$\fi}
\def\newcount{\alloc@0\count\countdef\insc@unt}
\def\newdimen{\alloc@1\dimen\dimendef\insc@unt}
\def\newwrite{\alloc@7\write\chardef\sixt@@n}

\newwrite\@unused
\newcount\@tempcnta
\newcount\@tempcntb
\newdimen\@tempdima
\newdimen\@tempdimb
\newbox\@tempboxa

\def\@spaces{\space\space\space\space}
\def\@whilenoop#1{}
\def\@whiledim#1\do #2{\ifdim #1\relax#2\@iwhiledim{#1\relax#2}\fi}
\def\@iwhiledim#1{\ifdim #1\let\@nextwhile=\@iwhiledim
        \else\let\@nextwhile=\@whilenoop\fi\@nextwhile{#1}}
\def\@badlinearg{\@latexerr{Bad \string\line\space or \string\vector
   \space argument}}
\def\@latexerr#1#2{\begingroup
\edef\@tempc{#2}\expandafter\errhelp\expandafter{\@tempc}%
\def\@eha{Your command was ignored.
^^JType \space I <command> <return> \space to replace it
  with another command,^^Jor \space <return> \space to continue without it.}
\def\@ehb{You've lost some text. \space \@ehc}
\def\@ehc{Try typing \space <return>
  \space to proceed.^^JIf that doesn't work, type \space X <return> \space
to
  quit.}
\def\@ehd{You're in trouble here.  \space\@ehc}

\typeout{LaTeX error. \space See LaTeX manual for explanation.^^J
 \space\@spaces\@spaces\@spaces Type \space H <return> \space for
 immediate help.}\errmessage{#1}\endgroup}
\def\typeout#1{{\let\protect\string\immediate\write\@unused{#1}}}

\font\tenln    = line10
\font\tenlnw   = linew10

\newdimen\@wholewidth
\newdimen\@halfwidth
\newdimen\unitlength

\unitlength =1pt


\def\thinlines{\let\@linefnt\tenln \let\@circlefnt\tencirc
  \@wholewidth\fontdimen8\tenln \@halfwidth .5\@wholewidth}
\def\thicklines{\let\@linefnt\tenlnw \let\@circlefnt\tencircw
  \@wholewidth\fontdimen8\tenlnw \@halfwidth .5\@wholewidth}

\def\linethickness#1{\@wholewidth #1\relax \@halfwidth .5\@wholewidth}

\newif\if@negarg

\def\lline(#1,#2)#3{\@xarg #1\relax \@yarg #2\relax
\@linelen=#3\unitlength
\ifnum\@xarg =0 \@vline
  \else \ifnum\@yarg =0 \@hline \else \@sline\fi
\fi}

\def\@sline{\ifnum\@xarg< 0 \@negargtrue \@xarg -\@xarg \@yyarg -\@yarg
  \else \@negargfalse \@yyarg \@yarg \fi
\ifnum \@yyarg >0 \@tempcnta\@yyarg \else \@tempcnta -\@yyarg \fi
\ifnum\@tempcnta>6 \@badlinearg\@tempcnta0 \fi
\setbox\@linechar\hbox{\@linefnt\@getlinechar(\@xarg,\@yyarg)}%
\ifnum \@yarg >0 \let\@upordown\raise \@clnht\z@
   \else\let\@upordown\lower \@clnht \ht\@linechar\fi
\@clnwd=\wd\@linechar
\if@negarg \hskip -\wd\@linechar \def\@tempa{\hskip -2\wd\@linechar}\else
     \let\@tempa\relax \fi
\@whiledim \@clnwd <\@linelen \do
  {\@upordown\@clnht\copy\@linechar
   \@tempa
   \advance\@clnht \ht\@linechar
   \advance\@clnwd \wd\@linechar}%
\advance\@clnht -\ht\@linechar
\advance\@clnwd -\wd\@linechar
\@tempdima\@linelen\advance\@tempdima -\@clnwd
\@tempdimb\@tempdima\advance\@tempdimb -\wd\@linechar
\if@negarg \hskip -\@tempdimb \else \hskip \@tempdimb \fi
\multiply\@tempdima \@m
\@tempcnta \@tempdima \@tempdima \wd\@linechar \divide\@tempcnta \@tempdima
\@tempdima \ht\@linechar \multiply\@tempdima \@tempcnta
\divide\@tempdima \@m
\advance\@clnht \@tempdima
\ifdim \@linelen <\wd\@linechar
   \hskip \wd\@linechar
  \else\@upordown\@clnht\copy\@linechar\fi}

\def\@hline{\ifnum \@xarg <0 \hskip -\@linelen \fi
\vrule height \@halfwidth depth \@halfwidth width \@linelen
\ifnum \@xarg <0 \hskip -\@linelen \fi}

\def\@getlinechar(#1,#2){\@tempcnta#1\relax\multiply\@tempcnta 8
\advance\@tempcnta -9 \ifnum #2>0 \advance\@tempcnta #2\relax\else
\advance\@tempcnta -#2\relax\advance\@tempcnta 64 \fi
\char\@tempcnta}

\def\vector(#1,#2)#3{\@xarg #1\relax \@yarg #2\relax
\@linelen=#3\unitlength
\ifnum\@xarg =0 \@vvector
  \else \ifnum\@yarg =0 \@hvector \else \@svector\fi
\fi}

\def\@hvector{\@hline\hbox to 0pt{\@linefnt
\ifnum \@xarg <0 \@getlarrow(1,0)\hss\else
    \hss\@getrarrow(1,0)\fi}}

\def\@vvector{\ifnum \@yarg <0 \@downvector \else \@upvector \fi}

\def\@svector{\@sline
\@tempcnta\@yarg \ifnum\@tempcnta <0 \@tempcnta=-\@tempcnta\fi
\ifnum\@tempcnta <5
  \hskip -\wd\@linechar
  \@upordown\@clnht \hbox{\@linefnt  \if@negarg
  \@getlarrow(\@xarg,\@yyarg) \else \@getrarrow(\@xarg,\@yyarg) \fi}%
\else\@badlinearg\fi}

\def\@getlarrow(#1,#2){\ifnum #2 =\z@ \@tempcnta='33\else
\@tempcnta=#1\relax\multiply\@tempcnta \sixt@@n \advance\@tempcnta
-9 \@tempcntb=#2\relax\multiply\@tempcntb \tw@
\ifnum \@tempcntb >0 \advance\@tempcnta \@tempcntb\relax
\else\advance\@tempcnta -\@tempcntb\advance\@tempcnta 64
\fi\fi\char\@tempcnta}

\def\@getrarrow(#1,#2){\@tempcntb=#2\relax
\ifnum\@tempcntb < 0 \@tempcntb=-\@tempcntb\relax\fi
\ifcase \@tempcntb\relax \@tempcnta='55 \or
\ifnum #1<3 \@tempcnta=#1\relax\multiply\@tempcnta
24 \advance\@tempcnta -6 \else \ifnum #1=3 \@tempcnta=49
\else\@tempcnta=58 \fi\fi\or
\ifnum #1<3 \@tempcnta=#1\relax\multiply\@tempcnta
24 \advance\@tempcnta -3 \else \@tempcnta=51\fi\or
\@tempcnta=#1\relax\multiply\@tempcnta
\sixt@@n \advance\@tempcnta -\tw@ \else
\@tempcnta=#1\relax\multiply\@tempcnta
\sixt@@n \advance\@tempcnta 7 \fi\ifnum #2<0 \advance\@tempcnta 64 \fi
\char\@tempcnta}

\def\@vline{\ifnum \@yarg <0 \@downline \else \@upline\fi}

\def\@upline{\hbox to \z@{\hskip -\@halfwidth \vrule
  width \@wholewidth height \@linelen depth \z@\hss}}

\def\@downline{\hbox to \z@{\hskip -\@halfwidth \vrule
  width \@wholewidth height \z@ depth \@linelen \hss}}

\def\@upvector{\@upline\setbox\@tempboxa\hbox{\@linefnt\char'66}\raise
     \@linelen \hbox to\z@{\lower \ht\@tempboxa\box\@tempboxa\hss}}

\def\@downvector{\@downline\lower \@linelen
      \hbox to \z@{\@linefnt\char'77\hss}}

\thinlines

\newcount\@xarg
\newcount\@yarg
\newcount\@yyarg
\newcount\@multicnt
\newdimen\@xdim
\newdimen\@ydim
\newbox\@linechar
\newdimen\@linelen
\newdimen\@clnwd
\newdimen\@clnht
\newdimen\@dashdim
\newbox\@dashbox
\newcount\@dashcnt
 \catcode`@=12


\newbox\tbox
\newbox\tboxa

\def\leftzer#1{\setbox\tbox=\hbox to 0pt{#1\hss}%
     \ht\tbox=0pt \dp\tbox=0pt \box\tbox}

\def\rightzer#1{\setbox\tbox=\hbox to 0pt{\hss #1}%
     \ht\tbox=0pt \dp\tbox=0pt \box\tbox}

\def\centerzer#1{\setbox\tbox=\hbox to 0pt{\hss #1\hss}%
     \ht\tbox=0pt \dp\tbox=0pt \box\tbox}

%
\def\image(#1,#2)#3{\vbox to #1{\offinterlineskip
    \vss #3 \vskip #2}}


\def\leftput(#1,#2)#3{\setbox\tboxa=\hbox{%
    \kern #1\unitlength
    \raise #2\unitlength\hbox{\leftzer{#3}}}%
    \ht\tboxa=0pt \wd\tboxa=0pt \dp\tboxa=0pt\box\tboxa}

\def\rightput(#1,#2)#3{\setbox\tboxa=\hbox{%
    \kern #1\unitlength
    \raise #2\unitlength\hbox{\rightzer{#3}}}%
    \ht\tboxa=0pt \wd\tboxa=0pt \dp\tboxa=0pt\box\tboxa}

\def\centerput(#1,#2)#3{\setbox\tboxa=\hbox{%
    \kern #1\unitlength
    \raise #2\unitlength\hbox{\centerzer{#3}}}%
    \ht\tboxa=0pt \wd\tboxa=0pt \dp\tboxa=0pt\box\tboxa}

\unitlength=1mm

\def\put(#1,#2)#3{\noalign{\nointerlineskip
                               \centerput(#1,#2){$#3$}
                                \nointerlineskip}}


\def\cal#1{\Cal#1} 
\def\bboard{\Bbb} 

\font\sevenbboard=msbm7
\font\tengoth=eufm10
\font\sevengoth=eufm7

\def\setA{{\bboard A}}
\def\setAA{\hbox{\sevenbboard A}}
\def\setZ{{\bboard Z}}
\def\calC{{\cal C}}
\def\calB{{\cal B}}
\def\calA{{\cal A}}
\def\calI{{\cal I}}
\def\Sym{\hbox{\tengoth S}}
\def\SSym{\hbox{\sevengoth S}}
\def\modI{\!\!\!\pmod{\calI}}
\def\modIq{\!\!\!\pmod{\calI_q}}
\def\Ferm{\mathop{\text{\rm Ferm}}\nolimits}
\def\Bos{\mathop{\text{\rm Bos}}\nolimits}

\def\len{{\ell}}
\def\inv{\mathop{\text{\rm inv}}\nolimits}
\def\imv{\mathop{\text{\rm imv}}\nolimits}
\def\bw#1#2{\bigl(\smallmatrix{#1\cr#2\cr}\bigr)}
\def\sbw#1#2{\bigl[\smallmatrix{#1\cr#2\cr}\bigr]}
\def\bww#1#2#3#4#5#6{\bigl[\sbw{#1#2}{#4#5} \bw{#3}{#6} \bigr]}
\def\wbw#1#2#3#4#5#6{\bigl[\bw{#1}{#4} \sbw{#2#3}{#5#6} \bigr]}

\def\smallmatrix#1{\vcenter{\offinterlineskip
     \halign{\vrule height 4pt depth 2pt width 0pt
      \hfil$\scriptstyle##$\hfil&&\kern
       3pt\hfil$\scriptstyle##$\hfil \crcr#1\crcr}}}

\catcode`\@=11
\def\matrice#1{\null \,\vcenter {\normalbaselines \m@th
\ialign {\hfil $##$\hfil &&\  \hfil $##$\hfil\crcr
\mathstrut \crcr \noalign {\kern -\baselineskip } #1\crcr
\mathstrut \crcr \noalign {\kern -\baselineskip }}}\,}

\def\pmatrice#1{\left(\null\vcenter {\normalbaselines \m@th
\ialign {\hfil $##$\hfil &&\thinspace  \hfil $##$\hfil\crcr
\mathstrut \crcr \noalign {\kern -\baselineskip } #1\crcr
\mathstrut \crcr \noalign {\kern -\baselineskip }}}\right)}

\def\bmatrice#1{\left[\null\vcenter {\normalbaselines \m@th
\ialign {\hfil $##$\hfil &&\thinspace  \hfil $##$\hfil\crcr
\mathstrut \crcr \noalign {\kern -\baselineskip } #1\crcr
\mathstrut \crcr \noalign {\kern -\baselineskip }}}\right]}

\catcode`\@=12

\frenchspacing

\def\RefBe{Be78}
\def\RefFH{FH05}
\def\RefGLZ{GLZ05}
\def\RefKa{Ka95}
\def\RefLT{LT06}
\def\RefLoI{Lo83}
\def\RefLoII{Lo02}
\def\RefMa{Ma15}
\def\RefPW{PW91}
\def\RefRTI{RT02}
\def\RefRTII{RT05}

\def\RefBe{1}
\def\RefFH{2}
\def\RefGLZ{3}
\def\RefKa{4}
\def\RefLT{5}
\def\RefLoI{6}
\def\RefLoII{7}
\def\RefMa{8}
\def\RefPW{9}
\def\RefRTI{10}
\def\RefRTII{11}



\topmatter
\title A basis for the right quantum algebra
and the ``$1=q$" principle
\endtitle

\author Dominique Foata and Guo-Niu Han\endauthor

\address Institut Lothaire, 1 rue Murner, F-67000 Strasbourg, France \endaddress
\email foata\@math.u-strasbg.fr \endemail

\address I.R.M.A. UMR 7501, Universit\'e Louis Pasteur et CNRS, 
7 rue Ren\'e-Descartes, F-67084 Strasbourg, France \endaddress
\email guoniu\@math.u-strasbg.fr \endemail


\date January 1, 1994 and, in revised form, June 22, 1994\enddate


\abstract
We construct a basis for the right quantum algebra introduced
by Garoufalidis, L\^e and Zeilberger and give a method
making it possible to go from an algebra submitted to
commutation relations (without the variable~$q$) to the
right quantum algebra by means of an appropriate
weight-function. As a consequence, a strong quantum
MacMahon Master Theorem is derived. Besides, the algebra of
biwords is systematically in use.
\endabstract

\endtopmatter

\document

\head 1. Introduction\endhead 
In their search for a {\it natural $q$-analogue} of the
MacMahon Master Theorem Garoufalidis {\it et al.} [\RefGLZ]
have introduced the {\it right quantum algebra}~${\cal R}_q$
defined to be the associative algebra over a commutative
ring~$K$, generated by $r^2$ elements
$X_{xa}$ $(1\le x,a\le r)$ ($r\ge 2$) subject to the following
commutation relations:
$$
\vbox{\hfil$\eqalign{
X_{yb} X_{xa} - X_{xa} X_{yb}
&= q^{-1}X_{xb}
X_{ya}-qX_{ya} X_{xb},\cr
{X}_{ya} X_{xa}&=q^{-1} X_{xa}
X_{ya},\cr}\quad
\eqalign{
&\hbox{($x>y$, $a>b$)}; \cr
&(x>y,\ \hbox{all }a); \cr
}$\hfil
}
$$
with $q$ belonging to~$K$.
The right quantum algebra in the case $r=2$ has already
been studied by Rodr\'\i guez-Romo and Taft [\RefRTI], 
who set up an explicit basis for it. On the other
hand, a basis for the {\it full quantum} algebra has been duly
constructed (see [\RefPW, Theorem 3.5.1, p. 38]) for an
arbitrary $r\ge 2$. It then seems natural to do the same with
the right quantum algebra for each $r\ge 2$. This is the
first goal of the paper.

In fact, the paper originated from a discussion with
Doron Zeilberger, when he explained to the first author how
he verified the quantum MacMahon Master identity for each
fixed~$r$ by computer code. His computer program uses the
fact, as he is perfectly aware, that the set of {\it irreducible
biwords}---which will be introduced in the sequel---{\it
generates} the right quantum algebra. For a better
understanding of his joint paper [\RefGLZ] and also for
deriving  the ``$1=q$" principle, it seems
essential to see whether the set of irreducible biwords has
the further property of being a {\it basis} and it does.
Thanks to this result, a {\it strong quantum
MacMahon Master Theorem} can be further derived.

\medskip
When manipulating the above relations, the visible part of the
commutations is made within the subcripts of the
$X_{xa}$'s, written as such. It is then important to magnify
them by having an adequate notation. To that end we replace
each
$$X_{xa}\hbox{ by the
{\it biletter} }\textstyle{x\choose a},
$$
so that, in further computations,
products $X_{x_1,a_1}X_{x_2,a_2}\cdots X_{x_\ell,a_\ell}$
become {\it biwords}
${x_1 x_2 \cdots x_\ell\choose a_1 a_2 \cdots a_\ell}$,
objects that have been efficiently used
in combinatorial contexts [\RefLoI, \RefLoII].
Accordingly, the commutation rules for the right
quantum algebra reread:
$$\eqalign{\textstyle
{yx\choose ba}-{xy\choose ab}&=
\textstyle  q^{-1}{xy\choose ba}-q{yx\choose ab}
\kern1.75cm (x>y,\, a>b);\cr
\textstyle{yx\choose aa}
&= \textstyle q^{-1}{xy\choose
aa}
\kern3cm (x>y,\; {\text{all}\ }a).\cr }
$$

We shall use the following notations.
The positive integer $r$ will be kept fixed throughout
and $\setA$ will designate the {\it alphabet}
$\{1,2,\ldots,r\}$. A {\it biword} on~$\setA$  is
a $2\times n$ matrix~$\alpha={x_1\cdots x_n\choose
a_1\cdots a_n}$ ($n\geq 0$), whose entries are in
$\setA$, the first (resp. second) row being called
the {\it top word} (resp. {\it bottom word\/}) of the
biword~$\alpha$. The number $n$ is the {\it length} of
$\alpha$; we write $\len(\alpha)=n$. The biword~$\alpha$
can also be viewed as a {\it word} of {\it biletters}
${x_1\choose a_1}\cdots{x_n\choose a_n}$, those biletters
$x_i\choose a_i$ being pairs of integers written vertically
with $x_i,a_i\in \setA$ for all $i=1,\ldots,n$. The {\it
product} of two biwords is their concatenation.

Let ${\bboard B}$ denote the set of biletters and $\calB$
the set of all biwords.
Let $\alpha\in \calB$ such that
$\len(\alpha)\geq 2$
and $1\leq i\leq
\len(\alpha)-1$ be a positive integer. The biword~$\alpha$
can be factorized as
$\alpha=\beta\bw{xy}{ab}\gamma$,  where
$x,y,a,b\in \setA$ and $\beta,\gamma\in \calB$ with
$\len(\beta)=i-1$. Say that $\alpha$ has a {\it double
descent} at position $i$ if $x>y$ and $a\geq b$. Notice the
discrepancy in the inequalities on the top word and the
bottom one. A biword~$\alpha$ without any double descent is
said to be {\it irreducible}. The set of all irreducible
biwords is denoted by ${\cal B}_{\text{irr}}$.

Let $\bboard Z$ be the ring of all integers and
${\bboard Z}[q,q^{-1}]$ the ring of the
polynomials in the variables~$q$, $q^{-1}$ submitted to
the rule $qq^{-1}=1$ with integral coefficients. The set
$\calA={\bboard Z}\langle\!\langle \calB\rangle\!\rangle$
of the formal sums $\sum_{\alpha} c(\alpha)
\alpha$, where $\alpha\in\calB$ and $c(\alpha)\in\setZ$,
together with the above biword multiplication,
the free addition and the free scalar product
forms an algebra over $\setZ$, called the {\it free biword
large $\setZ$-algebra}. The formal sums
$\sum_\alpha c(\alpha)\alpha$ will be called {\it
expressions}. An expression $\sum_\alpha c(\alpha)
\alpha$ is said to be  {\it irreducible} if $c(\alpha)=0$ for
all $\alpha\not\in\calB_{\text{irr}}$. The set of all
irreducible expressions is denoted by
${\cal A}_{\text{irr}}$.

Similarly, 
let $\calA_q={\bboard Z}[q,q^{-1}]
\langle\!\langle\calB\rangle\!\rangle$
denote the large ${\bboard Z}[q,q^{-1}]$-algebra of the
formal sums $\sum_\alpha c(\alpha)\alpha$, where
$c(\alpha)\in {\bboard Z}[q,q^{-1}]$ for all $\alpha\in {\cal B}$.
Following Bergman's method [\RefBe] we introduce two {\it
reduction systems}
$(S)$ and
$(S_q)$ as being the sets of pairs
$(\alpha,[\alpha])\in {\cal B}\times {\cal A}$ and
$(\alpha,[\alpha]_q)\in {\cal B}\times {\cal A}_q$,
respectively, defined by
$$
\eqalign{\textstyle
\bigl({xy\choose ab}, {xy\brack ab}\bigr)&\hbox{ with }
\textstyle{xy\brack ab}=  {yx\choose ba} + {yx\choose ab} -
{xy\choose ba}
\kern1.5cm (x>y,\, a>b);\cr
\textstyle\bigl({xy\choose aa},{xy\brack aa}\bigr)
&\hbox{ with }\textstyle {xy\brack aa}= {yx\choose aa}
\kern3.66cm (x>y,\; {\text{all}\ }a).\cr }\leqno(S)
$$
\vskip-5pt
$$
\eqalign{\textstyle
\bigl({xy\choose ab}, {xy\brack ab}_q\bigr)&\hbox{ with }
\textstyle{xy\brack ab}_q=  {yx\choose ba} + q{yx\choose
ab} - q^{-1}{xy\choose ba}
\quad (x>y,\, a>b);\cr
\textstyle\bigl({xy\choose aa},{xy\brack aa}_q\bigr)
&\hbox{ with }\textstyle {xy\brack aa}_q= q{yx\choose
aa}
\kern3cm (x>y,\; {\text{all}\ }a).\cr }\leqno(S_q)
$$
Notice that the equations ${xy\choose
ab}-{xy\brack ab}_q=0$, ${xy\choose aa}-{xy\brack aa}_q=0$
are simple rewritings of the commutation rules of the right
quantum algebra
and that
$(S)$ is deduced from
$(S_q)$ by letting $q=1$.

Let 
$\calI$ (resp.~$\calI_q$) be the two-sided ideal
of~$\calA$ (resp. of~$\calA_q$) generated by the elements
$\gamma-[\gamma]$ (resp. $\gamma-[\gamma]_q$) such
that $(\gamma,[\gamma])\in (S)$ (resp.
$\gamma-[\gamma]_q\in (S_q)$). The quotient algebras
${\cal R}={\cal A}/\calI$ and ${\cal R}_q={\cal A}_q/\calI_q$ are
called the {\it $1$-quantum right algebra} and the {\it
$q$-right quantum algebra}, respectively.

In Section~2 we define a $\setZ$-linear mapping $E\mapsto [E]$
of~$\cal A$ onto the ${\bboard Z}$-module~${\cal A}_{\text
irr}$ of the irreducible expressions, called itself {\it
reduction}. Using Bergman's ``Diamond Lemma" [\RefBe] this
reduction will serve to obtain a model for the algebra~$\cal
R$, as stated in the next theorem.

\goodbreak
\proclaim{Theorem 1} A set of representatives in~$\cal A$ 
for~$\cal R$ is given
by the ${\bboard Z}$-module~${\cal A}_{\text{irr}}$. The
algebra~$\cal R$ may be identified with the ${\bboard
Z}$-module~${\cal A}_{\text{irr}}$, the multiplication being
given by $E\times F=[E\times F]$ for any two irreducible
expressions $E$, $F$. With this identification ${\cal
B}_{\text{irr}}$ is a basis for~${\cal R}$.
\endproclaim

Four statistics counting various kinds of {\it inversions} will
be needed. If $\alpha=\bw uv=\bw{x_1x_2\ldots
x_n}{a_1a_2\ldots a_n}$ is a biword, let
$$
\eqalignno{
\inv u &=\#\{(i,j) \mid 1\leq i< j\leq n, x_i>x_j \}; \cr
\imv v &=\#\{(i,j) \mid 1\leq i< j\leq n, a_i\geq a_j \}; \cr
\inv^- (\alpha) &=\inv v-\inv u; \cr
\inv^+ (\alpha) &=\imv v+\inv u. \cr
}
$$
The first (resp. second) statistic ``inv" (resp. ``imv") is the usual number
of inversions (resp. of {\it large} inversions) of a word.
Notice that ``$\inv^-$" may be negative.
The weight function~$\phi$ defined for each biword~$\alpha$
by $\phi(\alpha)=q^{\inv^-(\alpha)}\alpha$ can be extended
to all of~${\cal A}_q$ by linearity.
Clearly it is
a  $\setZ[q,q^{-1}]$-module isomorphism
of~${\cal A}_q$ onto itself. The second result of the paper
is stated next.

\proclaim{Theorem 2} The
$\setZ[q,q^{-1}]$-module isomorphism of $\calA_q$ onto
itself
$\phi$ induces a
${\bboard Z}[q,q^{-1}]$-module isomorphism
$\overline \phi$ of
$\calA_q/\calI$
onto $\calA_q/\calI_q$. In particular,~${\cal R}_q$ has the same
basis as~$\cal R$.
\endproclaim

Now, a {\it circuit} is defined to be a biword whose top word
is a rearrangement of the letters of its bottom word.  The
set of all circuits is denoted by~$\calC$. It is clearly a
submonoid of~$\calB$. An expression
$E=\sum_\alpha c(\alpha)\alpha$ is said to be {\it circular}
if  $c(\alpha) =0$ for all $\alpha\not\in\calC$. Clearly, the
sum and the product of two circular expressions is a circular
expression, so that the set of circular expressions in~${\cal
A}$ (resp. in~${\cal A}_q$) is a subalgebra of~${\cal
A}$ (resp. of~${\cal A}_q$), which
will be denoted by 
$\calA^{\text{cir}}$ (resp. by~$\calA_q^{\text{cir}}$).

\proclaim{Theorem 3} The restriction of the 
${\bboard Z}[q,q^{-1}]$-module
isomorphism~$\overline\phi$ of Theorem~$2$ to
$\calA_q^{\text{cir}}/\calI$ is a
${\bboard Z}[q,q^{-1}]$-algebra isomorphism onto
$\calA_q^{\text{cir}}/\calI_q$.
\endproclaim
\vskip.1cm

\def\fleche(#1,#2)\dir(#3,#4)\long#5{%
{\leftput(#1,#2){\vector(#3,#4){#5}}}}

$$\hbox{\vbox{\offinterlineskip
\def\Acir{{\cal A}^{\text{cir}}_q}
\medskip
\centerput(0,0){$\Acir$}
\fleche(8,0)\dir(1,0)\long{17}
\centerput(32,0){$\Acir$}
\centerput(16,2){$\phi$}
\centerput(16,-3){\sevenrm module-iso}
\centerput(-3,-8){\sevenrm can}
\fleche(0,-3)\dir(0,-1)\long{11}
\centerput(0,-18){$\Acir/{\cal I}$}
\fleche(8,-17)\dir(1,0)\long{17}
\centerput(32,-18){$\Acir/{\cal I}_q$}
\centerput(16,-15){$\overline \phi$}
\centerput(16,-20){\sevenrm algebra-iso}
\fleche(32,-3)\dir(0,-1)\long{11}
\centerput(29,-8){\sevenrm can}
\vskip 1.5cm
}\hskip6cm
\vbox{\offinterlineskip
\medskip
\centerput(0,0){${\cal A}_q$}
\fleche(8,0)\dir(1,0)\long{17}
\centerput(32,0){${\cal A}_q$}
\centerput(16,2){$\phi$}
\centerput(16,-3){\sevenrm module-iso}
\centerput(-3,-8){\sevenrm can}
\fleche(0,-3)\dir(0,-1)\long{11}
\centerput(0,-18){${\cal A}_q/{\cal I}$}
\fleche(8,-17)\dir(1,0)\long{17}
\centerput(32,-18){${\cal A}_q/{\cal I}_q$}
\centerput(16,-15){$\overline \phi$}
\centerput(16,-20){\sevenrm module-iso}
\fleche(32,-3)\dir(0,-1)\long{11}
\centerput(29,-8){\sevenrm can}
\centerput(-15,-8){$\hookrightarrow$}
\vskip 1.5cm
}}\hskip3cm
$$
\vskip 0.7cm

Accordingly, each identity holding in
$\calA_q^{\text{cir}}/\calI$
has an equivalent counterpart in
$\calA_q^{\text{cir}}/\calI_q$.
This is the
``$1=q$" principle. An illustration of this principle is given
with the $q$MM Theorem derived by Garoufalidis et al. ({\it
op. cit.}). Those authors have introduced the
$q$-Fermion and the
$q$-Boson as being the sums
$$\eqalignno{
\Ferm(q)
&=\sum_{J\subset \setAA}
(-1)^{|J|} \sum_{\sigma\in {\SSym}_J}
(-q)^{-\inv \sigma}
\pmatrice  {\sigma (i_1)&\sigma (i_2)&\cdots&\sigma (i_l)\cr
i_1&i_2&\cdots&i_l\cr}\cr
\noalign{\hbox{and}}
\Bos(q)&=\sum_w  q^{\inv w}\pmatrice {
\overline w\cr w\cr},\cr}
$$
where ${\Sym}_J$ is the
permutation group acting on the set~$J=\{i_1<i_2<\cdots <i_l\}$
and where $\overline w$ stands for the nondecreasing rearrangement of
the word~$w$. The $q$-Fermion and $q$-Boson belong
to~${\cal A}_q$. Garoufalidis {\it et al.} have
proved ([\RefGLZ],  see also [\RefFH] for another proof)
the following 
identity
$$
\Ferm(q)\times \Bos(q) \equiv 1 \!\!\! \pmod{\calI_q}.
$$
This identity  may {\it a priori\/} be  regarded as
a {\it $q$-version of the MacMahon Master Theorem}, as it
reduces to the classical MacMahon's identity
[\RefMa, p.~93--98], when
$q=1$ and the biletters are supposed to commute.
By the ``$1=q$" principle (Theorem 3), we obtain
the following result, which shows that the variable~$q$ is in fact
superfluous.

\proclaim{Corollary 4} The identity
$\Ferm(q)\times \Bos(q)\equiv 1\modIq$
holds
if and only if
$\Ferm(1)\times \Bos(1)\equiv 1\modI$
holds.
\endproclaim

When applying Theorem 1 to the quantum MacMahon Master
Theorem, we obtain the following result, which can be
regarded as a {\it strong quantum MacMahon Master Theorem}.

\proclaim{Corollary 5}
The following identity holds:
$$[\Ferm(1)\times \Bos(1)] = 1.$$
\endproclaim

The proof of Theorem~1, given in section~2, is based on
Bergman's ``Diamond Lemma" ([\RefBe], Theorem~1.2). It
consists of verifying that the conditions required by the
``Diamond Lemma" hold in the present situation. Theorem~2
and~3 are proved in section~3 and Corollaries 4 and~5
in the last section.

\head 2. The Reduction Process\endhead 

\def\Birr{{\cal B}_{\text{irr}}}

For applying Bergman's ``Diamond Lemma" [\RefBe] we need
have a so-called {\it reduction system} that satisfies two
conditions: (1) the {\it descending chain condition} holds;
(2) its {\it ambiguities are resolvable}. Those terms will be
explained shortly. Let ${\bboard Z}\langle{\cal B}\rangle$ be
the subalgebra of~$\cal A$ of all
sums $\sum_\alpha c(\alpha)\alpha$ $(\alpha\in {\cal B})$,
where only a finite number of the $c(\alpha)$'s are
non-zero. The elements of ${\bboard Z}\langle{\cal
B}\rangle$ will be called {\it finite} expressions.
The reduction system for ${\bboard
Z}\langle{\cal B}\rangle$ is simply the set~$(S)$ of pairs
$(\alpha,[\alpha])\in {\cal B}\times {\bboard
Z}\langle{\cal B}\rangle$, already
described in section~1. By means of this reduction system we
introduce the {\it reduction} itself as being the mapping
$E\mapsto [E]$ of
${\bboard
Z}\langle{\cal B}\rangle$ onto the set ${\bboard
Z}\langle\Birr\rangle$ of the {\it finite}
irreducible expressions, defined by the following axioms:

(C1) The reduction is linear, i.e. for $E_1, E_2\in{\bboard
Z}\langle{\cal B}\rangle$,
$$[c_1 E_1+c _2 E_2]=c_1 [E_1]+c_ 2[E_2].$$

(C2) 
$[E]=E$ for every $E\in{\bboard
Z}\langle\Birr\rangle$.
\smallskip

(C3) For $\bw{xy}{ab}\not\in\Birr$, i.e. $x>y$ and $a\geq
b$, then
$$
\leqalignno{
\textstyle  {xy\brack ab} &=
\textstyle  {yx\choose ba} + {yx\choose ab} - {xy\choose
ba}, \hbox{\quad \rm if
$a>b$;}
&\hbox{\qquad\rm (C3.1)}\cr
  \textstyle{xy\brack aa} &=
 \textstyle {yx\choose aa},\hfill \hbox{\quad\rm if $a=b$.}
&\hbox{\qquad\rm (C3.2)}\cr
}
$$

(C4) Let $\alpha\not\in \Birr$ be a reducible biword.
For each factorization
$\alpha=\beta\bw{xy}{ab}\gamma$, where
$x,y,a,b\in \setA$ and $\beta\bw{x}{a}\in\Birr$
and $\gamma\in \calB$ we have
$$
[\alpha]=\bigl[\beta\times \sbw{xy}{ab}\times \gamma\bigr].
$$

Although the reduction is recursively defined, it is
well defined. Every time condition (C3) is applied,
the running biword $\alpha$ is transformed into either three
new biwords (C3.1), or a new biword (C3.2). The important
property is the fact that the statistic $\inv^+$ of each of
the new biwords is {\it strictly less} than $\inv^+(\alpha)$.
Thus, after {\it finitely} many successive applications of (C3),
an irreducible expression is derived. When the
biword~$\alpha$ has more than one double descent,
condition (C4) says that condition (C3) must be applied at
the first double descent position. Consequently, the final
irreducible expression is unique.

\goodbreak
There are several other ways to map each expression
onto an irreducible expression, using relations (C3.1)
and~(C3.2). The reduction defined above is only one of those
mappings. The important feature is the fact such a mapping
involves {\it finitely} many applications of relations~(C3.1)
and~(C3.2). Following Bergman we say that condition~(1)
({\it the descending chain condition}) holds for the pair
$({\bboard
Z}\langle{\cal B}\rangle,S)$.

Now an expression~$E$ is said to be
{\it reduction-unique under~$S$} if the irreducible
expression derived from~$E$ does not depend on {\it where}
the applications of (C3.1) and (C3.2) take place. More formally,
each expression  is reduction-unique
under~$S$ if for every biword~$\alpha$ and every
factorization
$\alpha=\beta\alpha'\gamma$ $(\beta, \alpha',\gamma\in
{\cal B})$ we have
$$
[\alpha]=[\beta[\alpha']\gamma].\kern2cm
\leqno\hbox{(Reduction-unique)}
$$

Now examine the second condition required by the ``Diamond
Lemma." There is an {\it
ambiguity} in~$S$ if there are two pairs
$(\alpha,[\alpha])$ and
$(\alpha',[\alpha'])$ in~$S$ such that
$\alpha=\beta\gamma$,
$\alpha'=\gamma\delta$ for some nonempty biwords
$\beta,\gamma,\delta\in {\cal B}$. This ambiguity is said to
be {\it resolvable} if $[[\alpha]\delta]=[\beta[\alpha']]$, i.e.
$[[\beta\gamma]\delta]=[\beta[\gamma\delta]]$. In our
case this means that $\beta={x\choose a}$,
$\gamma={y\choose b}$, $\delta={z\choose c}$ with $x>y>z$
and $a\ge b\ge c$. Using the first three integers instead of
the letters $x,y,z,a,b,c$ there are four cases to study:

(i) $\beta={3\choose 3}$, $\gamma={2\choose 2}$,
$\delta ={1\choose 1}$;\qquad

(ii) $\beta={3\choose 2}$, $\gamma={2\choose 2}$,
$\delta ={1\choose 1}$;\qquad

(iii) $\beta={3\choose 2}$, $\gamma={2\choose 1}$,
$\delta ={1\choose 1}$;\quad

(iv) $\beta={3\choose 1}$, $\gamma={2\choose 1}$,
$\delta ={1\choose 1}$.
\bigskip

Accordingly, the ambiguities in~$S$ are {\it resolvable} if
the following four identities hold:
$$\leqalignno{
\textstyle [{32\brack 32}{1\choose 1}]
&=\textstyle[{3\choose 3}{21\brack 21}];&(2.1)\cr
\textstyle [{32\brack 22}{1\choose 1}]
&=\textstyle[{3\choose 2}{21\brack 21}];&(2.2)\cr
\textstyle [{32\brack 21}{1\choose 1}]
&=\textstyle[{3\choose 2}{21\brack 11}];&(2.3)\cr
\textstyle [{32\brack 11}{1\choose 1}]
&=\textstyle[{3\choose 1}{21\brack 11}].&(2.4)\cr
}
$$
Let us prove those four identities.

\medskip
{\it Proof} of (2.1).\quad
For an easy reading of the coming calculations we have added
subscripts $A$, $B$, \dots~, $J$ to certain brackets, which
should help spot the subscripted brackets in the various
equations. The left-hand side is evaluated as follows:
$$
\eqalignno{  
\bww321321   &= \bww231231_A
+\bww231321_B -\bww321231_C;\cr
\bww231231_A &= \wbw231231
= \wbw213213_D +\wbw213231_I -\wbw231213_J;\cr
\bww231321_B &= \wbw231321
= \wbw213312_E +\wbw213321_F -\wbw231312;\cr
-\bww321231_C &= -\wbw321231
= -\wbw312213_G -\wbw312231 +\wbw321213_H;\cr
\wbw213213_D &= \bww213213
= \bww123123 +\bww123213
-\bww213123;\cr}
$$

$$\eqalignno{  
\wbw213312_E &= \bww213312
= \bww123132 +\bww123312 -\bww213132;\cr
\wbw213321_F &= \bww213321
= \bww123231 +\bww123321 -\bww213231_I;\cr
-\wbw312213_G &= -\bww312213
= -\bww132123 -\bww132213 +\bww312123;\cr
\wbw321213_H &= \bww321213
= \bww231123 +\bww231213_J
-\bww321123.\cr }$$
The sum of the above nine equalities yields:
$$
\eqalign{\bww321321=
&-\bw{231}{312}-\bw{312}{231}+\bw{123}{123}
+\bw{123}{213}-\bw{213}{123}\cr
&+\bw{123}{132}+\bw{123}{312}-\bw{213}{132}
+\bw{123}{231}+\bw{123}{321}\cr
&-\bw{132}{123}-\bw{132}{213}+\bw{312}{123}
+\bw{231}{123}-\bw{321}{123}.\cr
}$$
As for the right-hand side we have:
$$
\eqalign{  
\wbw321321   &= \wbw312312_A +\wbw312321_B
-\wbw321312_C;\cr
\bww312312_A &= \bww132132_D + \bww132312_I
- \bww312132_J ;\cr
\bww312321_B &= \bww132231_E + \bww132321_F
- \bww312231 ;\cr
-\bww321312_C &= -\bww231132_G- \bww231312
+ \bww321132_H ;\cr
\wbw132132_D &= \wbw123123 + \wbw123132
-\wbw132123 ;\cr
\wbw132231_E &= \wbw123213 + \wbw123231
- \wbw132213   ;\cr
\wbw132321_F &=\wbw123312  +\wbw123321
-\wbw132312_I  ;\cr
-\wbw231132_G &=-\wbw213123  -\wbw213132
+\wbw231123  ;\cr 
\wbw321132_H &=\wbw312123  +\wbw312132_J
-\wbw321123.  \cr 
}$$
The sum of the above nine equalities yields:
$$
\eqalign{  
\wbw321321 = 
&-\bw{231}{312}-\bw{312}{231}+\bw{123}{123}+\bw{123}{213}-\bw{213}{123}\cr
&+\bw{123}{132}+\bw{123}{312}-\bw{213}{132}+\bw{123}{231}+\bw{123}{321}\cr
&-\bw{132}{123}-\bw{132}{213}+\bw{312}{123}+\bw{231}{123}-\bw{321}{123}\cr
&= \bww321321 \ \hbox{(left-hand side)}.\cr
}$$

{\it Proof} of (2.2).\quad We have:
$$
\eqalignno{
\bww321221 &= \wbw231221= \wbw213212 +\wbw213221
-\wbw231212 \cr &= \bww123122 +\bww123212
-\bww213122 +\wbw213221 -\wbw231212 \cr &=
\bw{123}{122} +\bw{123}{212} -\bw{213}{122}
+\bw{123}{221} -\bw{231}{212} \cr
\noalign{\hbox{On the other hand,}}
\wbw321221 &=\wbw312212_A+\wbw312221_B-\wbw321212_C\cr
\bww312212_A &=\bww132122  +\bww132212_D  -\bww312122_E  \cr
\bww312221_B &=\wbw132221 =   \wbw123212  +\wbw123221  -\wbw132212_D  \cr
-\bww321212_C &=-\bww231122  -\bww231212  +\bww321122_E  \cr}
$$

\noindent
The sum of the above four equalities yields:
$$
\eqalignno{
\wbw321221 &= \bww132122 + \wbw123212  +\wbw123221
-\bww231122  -\bww231212\cr
&  =\bww321221\ \hbox{(left-hand side)}.\cr}
$$

{\it Proof} of (2.3).\quad We form
$$
\eqalignno{
\bww321211 &= \bww231121_A +\bww231211_B
-\bww321121_C ;\cr
\wbw231121_A &=\wbw213112  +\wbw213121_D
-\wbw231112_E  ;\cr
\wbw231211_B &=\bww213211 = \bww123121
+\bww123211  -\bww213121_D  ;\cr
-\wbw321121_C &=-\wbw312112  -\wbw312121
+\wbw321112_E; \cr 
\noalign{\hbox{so that}}
\bww321211  & =  
\bw{123}{112}  + \bw{123}{121}  +\bw{123}{211}
-\bw{132}{112}  -\bw{312}{121}.\cr
\noalign{\hbox{On the other hand,}}
\wbw321211&=\bww312211=\bww132121+\bww132211-\bww312121\cr
&=\wbw123112+\wbw123121-\wbw132112 + \bww132211-\bww312121 \cr
&= \bww321211 \ \hbox{(left-hand side)}. \cr
}$$

Finally, the proof of (2.4)
$[{32\brack 11}{1\choose 1}]
=[{3\choose 1}{21\brack 11}]$ is straightforward.\cqfd

\medskip
Let ${\cal I}'$ be the two-sided ideal of ${\bboard
Z}\langle{\cal B}\rangle$ generated by the elements
$\gamma-[\gamma]$ such that $(\gamma,[\gamma])\in
(S)$ and let ${\cal R}'$ be the quotient ${\bboard
Z}\langle{\cal B}\rangle/{\cal I}'$. The pair $({\bboard Z}\langle{\cal
B}\rangle,S)$ having the
descending chain condition and the ambiguities being
resolvable, Bergman's ``Diamond Lemma" implies the following
theorem.

\proclaim{Theorem 6}
\hskip -4pt A set of representatives in~${\bboard
Z}\langle{\cal B}\rangle$ for~${\cal R}'$ is given
by the ${\bboard Z}$-module~${\bboard Z}\langle\Birr\!\rangle\!$. 
The algebra~${\cal R}'$ may be identified with~${\bboard
Z}\langle\Birr\rangle$, the multiplication being
given by $E\times F=[E\times F]$ for any two
finite irreducible expressions $E$, $F$. With this
identification ${\cal B}_{\text{irr}}$ is a basis for~${\cal
R}'$.
\endproclaim

For the proof of Theorem~1 we use the following argument.
For each $n\ge 0$ let ${\cal A}^{(n)}$ be the $n$-th {\it
degree homogeneous subspace} of~$\cal A$ consisting of all
expressions $\sum_\alpha c(\alpha)\alpha$ such that
$\ell(\alpha)=0$ if $\ell(\alpha)\not=n$ ([\RefKa], I.6).
As the starting
alphabet~$\bboard A$ is finite, all expressions
in~${\cal A}^{(n)}$ are {\it finite} sums, so that
${\cal A}^{(n)}\subset {\bboard Z}\langle{\cal B}\rangle$ for
all~$n$. Now, let each expression~$E$ from~$\cal A$ be
written as a sum
$$
E=\sum_{n\ge 0} E^{(n)}\quad
\hbox{with } E^{(n)}\in {\cal A}^{(n)}.
$$
As both ideals $\cal I$ and ${\cal I}'$ are generated by
expressions from the second degree homogeneous
space~${\cal A}^{(2)}$, we have
$E\equiv 0\!\!\!\pmod{{\cal I}}$ if and only if for every
$n\ge 0$ we have $E^{(n)}\equiv 0\!\!\!\pmod{{\cal I}'}$.
In particular,
$[E^{(n)}]$ also belongs to ${\cal A}^{(n)}$. Theorem~1
follows from Theorem~6 by taking the definition
$$
[E]=\sum_{n\ge 0}[ E^{(n)}].
$$


\head 3. The ``$1=q$" principle\endhead 
Let $E$ be an expression in $\calA_q$.
Theorem 2 is equivalent to saying that
the identity $E\in\calI$ holds
if and only if $\phi(E) \in\calI_q$ holds.
First we prove the ``only if" part.
Because $\phi$ is linear, it suffices
to consider all linear generators of $\calI$, which have the following
form:
$$\displaylines{\noalign{\vskip-5pt}
\textstyle E_1=\alpha{rs\choose i\,i}\beta -  \alpha{s\,r\choose
i\,i}\beta\cr
\noalign{\vskip-5pt}
\noalign{\hbox{and}}
\textstyle E_2=\alpha{rs\choose ij}\beta - \alpha{sr\choose ji} \beta-
 \alpha{sr\choose ij}\beta + \alpha{rs\choose ji}\beta,\cr}$$
where $\alpha, \beta$ are biwords and $r$, $s$, $i$, $j$
integers such that $r<s$, $i<j$.
Let $k=\inv^-(\alpha{rs\choose ij}\beta)$,
then
$$\displaylines{\noalign{\vskip -5pt}
\textstyle \inv^-(\alpha{sr\choose ji}\beta)=k,\quad
\inv^-(\alpha{sr\choose ij}\beta)=k-1,\quad
\inv^-(\alpha{rs\choose ji}\beta)=k+1,\cr
\noalign{\hbox{so that}}
\textstyle \phi(E_2)=q^k\Bigl(
\alpha{rs\choose ij}\beta - \alpha{sr\choose ji} \beta-
 q^{-1}\alpha{sr\choose ij}\beta +q \alpha{rs\choose ji}\beta
\Bigr)\in\calI_q.\cr}
$$
In the same way $\phi(E_1)\in\calI_q$.
The ``if" part can be proved in the same manner.\qed
\medskip

For Theorem 3 it suffices to prove the identity
$\phi(EF)=\phi(E)\phi(F)$
for any two circular expressions $E$ and $F$. As $\phi$ is
linear, it suffices to do it when $E=\bw uv$,
$F=\bw {u'}{v'}$ are two circuits. Let $\inv(u,u')$ denote
the number of pairs
$(x,y)$ such that $x$ (resp.~$y$) is a letter of~$u$ (resp.
of~$u'$) and $x>y$. As $\bw {uu'}{vv'}$ is a circuit, we have
$$
\eqalignno{
\inv uu'&=\inv u + \inv u' + \inv(u,u') ; \cr
\inv vv'&=\inv v + \inv v' + \inv(v,v') ; \cr
\inv(u,u')&= \inv(v,v') ; \cr
\noalign{\hbox{so that}}
\inv^- (EF) &= \inv vv' -\inv uu' \cr
&= \inv v + \inv v' -\inv u - \inv u'\cr
&=\inv^-(E)+\inv^-(F). \qed\cr
}
$$

\head 4. The quantum MacMahon Master Theorem\endhead 

For proving Corollary 4 we apply Theorem 3
to $\Ferm(1)\times\Bos(1)$.
As $\Ferm(1)$ and $\Bos(1)$ are both circular expressions, the relation
$$
\displaylines{
\Ferm(1)\times \Bos(1)\equiv 1\modI\cr
\noalign{\hbox{is equivalent to}}
\phi(\Ferm(1))\times \phi(\Bos(1))\equiv 1\modIq.\cr
}$$
Finally, it is straightforward to verify
$$\Ferm(q)=\phi(\Ferm(1))
\quad{\text{and}}\quad\Bos(q)=\phi(\Bos(1)).
\qed$$ 

\medskip
Now to prove Corollary~5 we start with Garoufalidis {\it et
al.}'s result (for $q=1$), which says that
$\Ferm(1)\times \Bos(1)\equiv 1\modI$.
But by Theorem~1 we have
$E\equiv F \modI$ if and only if $[E]=[F]$. Hence
$[\Ferm(1)\times \Bos(1)]= [1]=1$.
\qed

\head 5. Concluding Remarks\endhead 
It is worth noticing that Rodr\'\i guez and Taft have introduced two explicit
left quantum groups for $r=2$ in [\RefRTI] and [\RefRTII]. As mentioned in the
introduction, the former one has been extended to an arbitrary dimension
$r\geq 2$ by Garoufalidis {\it et al.} [\RefGLZ] and been given an explicit 
basis in the present paper, while he latter one has been recently {\it
modeled after} ${\text{SL}}_q(r)$ by Lauve and Taft [\RefLT] also for
each $r\geq 2$.

We should like to thank Doron Zeilberger and Stavros Garoufalidis for
several fruitful discussions, orally and also by email. We are
grateful to Christian Kassel, who gave us the references to
the book by Parshall and Wang [\RefPW], the two papers by
Rodr\'\i guez and Taft [\RefRTI, \RefRTII] and the fundamental
paper by Bergman [\RefBe]. Finally, we have greatly benefited
from Jean-Pierre Jouanolou's ever-lasting mathematical
expertise. 


\refstyle{C}

\enddocument